\documentclass[a4paper, twoside, 11pt]{amsart}

\usepackage{geometry}                 
\usepackage{amsmath,amsfonts,amssymb} 
\usepackage{amsthm}                   
\usepackage{graphicx}                                  
\usepackage{dsfont}
\usepackage{enumerate}
\usepackage{hyperref}
\usepackage{enumitem}

\newtheorem{theorem}{Theorem}
\newtheorem*{theorem*}{Theorem}
\newtheorem{lemma}[theorem]{Lemma}

\newtheorem{corollary}[theorem]{Corollary}

\newcommand{\bz}{\mathbb{Z}}
\newcommand{\bt}{\mathbb{T}}
\newcommand{\br}{\mathbb{R}}
\newcommand{\norm}[1]{\left\lVert#1\right\rVert}
\newcommand{\1}{\mathds{1}}
\newcommand\dotprod[2]{\langle #1 , #2 \rangle}
\renewcommand{\geq}{\geqslant}
\renewcommand{\leq}{\leqslant}

\DeclareMathOperator{\mes}{mes}

\begin{document}

\pagenumbering{arabic}
\title{Multi-tiling and Riesz bases}
\author{Sigrid Grepstad}
\address{Department of Mathematical Sciences, Norwegian University of Science and Technology (NTNU), NO-7491 Trondheim, Norway.}
\email{\href{mailto:sigrid.grepstad@math.ntnu.no}{\texttt{sigrid.grepstad@math.ntnu.no}}}
\author{Nir Lev}
\address{Centre de Recerca Matem\`atica
Campus de Bellaterra, Edifici C, 
08193 Bellaterra (Barcelona), Spain}
\address{Department of Mathematics, Bar-Ilan University, Ramat-Gan 52900, Israel.}
\email{\href{mailto:levnir@math.biu.ac.il}{\texttt{levnir@math.biu.ac.il}}}
\thanks{Supported by the Israel Science Foundation grant No. 225/13}

\keywords{Riesz bases, Tiling, Quasicrystals}

\begin{abstract}
Let $S$ be a bounded, Riemann measurable set in $\br^d$, and $\Lambda$ be a lattice. By a theorem of Fuglede, if $S$ tiles $\br^d$ with translation set $\Lambda$, then $S$ has an orthogonal basis of exponentials. We show that, under the more general condition that $S$ \emph{multi-tiles} $\br^d$ with translation set $\Lambda$, $S$ has a Riesz basis of exponentials. The proof is based on Meyer's quasicrystals. 
\end{abstract} 

\maketitle
\section{Introduction}

\subsection{}
Let $S$ be a bounded, measurable set in $\br^d$. A sequence $\Lambda \subset \br^d$ is called a \textit{spectrum} for $S$ if the system of exponential functions
\begin{equation*}
 E(\Lambda) = \{e_\lambda\}_{\lambda \in \Lambda}, \quad \; e_\lambda(x) = e^{2 \pi i \dotprod{\lambda}{x}},
\end{equation*}
is a complete orthogonal system in $L^2(S)$. The existence of a spectrum $\Lambda$ for $S$ provides a unique and stable expansion of any $f \in L^2(S)$ into a ``non-harmonic'' Fourier series with frequencies in $\Lambda$.

The study of spectra was initiated by Fuglede \cite{fuglede} who suggested a connection with the concept of tiling. We say that $S$ tiles $\br^d$ with translation set $\Lambda$ if the sets $S + \lambda$ $(\lambda \in \Lambda)$ are disjoint and cover the whole space up to measure zero, that is,
\[
\sum_{\lambda \in \Lambda} \1_S(x - \lambda) = 1 \quad \text{(a.e.)}
\]

\begin{theorem*}[Fuglede \cite{fuglede}]
Let $\Lambda$ be a lattice. If $S$ tiles $\br^d$ with translation set $\Lambda$, then the dual lattice $\Lambda^*$ is a spectrum for $S$, and also the converse is true.
\end{theorem*}

By a lattice we mean the image of $\bz^d$ under some invertible linear transformation. The dual lattice $\Lambda^*$ is the set of all vectors $\lambda^* \in \br^d$ such that $\dotprod{\lambda}{\lambda^*} \in \bz$, $\lambda \in \Lambda$.

Fuglede conjectured that $S$ admits a spectrum if and only if it can tile $\br^d$ with some translation set. This conjecture inspired extensive research, see the surveys \cite{koloun3, koloun2} and the references therein for more information, including the present state of the conjecture.

\subsection{}
Fuglede's conjecture suggests that the existence of a spectrum is a rather special property of the set $S$. Indeed, it is known that some very simple sets $S$ do not admit an orthogonal basis of exponentials. For example, it is easy to construct a union of two intervals on $\br$ which does not admit such a basis (see also \cite{laba}). In dimension two, a convex polygon has a spectrum only if it is (up to an affine transformation) either a square or a hexagon \cite{iosevich}. The ball in any dimension $d \geq 2$ has no spectrum \cite{fuglede, iosevichped}.

The exponential system $E(\Lambda)$ is said to be a \emph{Riesz basis} for $S$ if the mapping $f \mapsto \{\dotprod{f}{e_\lambda}\}_{\lambda \in \Lambda}$ is bounded and invertible from $L^2(S)$ onto $\ell^2(\Lambda)$ (but not necessarily unitary, as in the case when $\Lambda$ is a spectrum). This condition still allows a unique and stable expansion of any $f \in L^2(S)$ into a Fourier series with frequencies in $\Lambda$, see \cite{young}, but it is not as rigid as the orthogonality requirement for the exponentials (for example, it is stable under small perturbations of the frequencies). Thus, one might hope that sets with no spectrum should at least have a Riesz basis of exponentials.

This was indeed established in some particular cases. It is a recent result due to Kozma and Nitzan \cite{kozmanitzan} that any finite union of intervals admits a Riesz basis of exponentials. In dimension two it was proved by Lyubarskii and Rashkovskii \cite{yurarash} that any convex, symmetric polygon has such a basis. 

In spite of this progress, however, there are still relatively few results on the existence of a Riesz basis of exponentials. In particular, it is still an open problem if such a basis exists for the ball in dimension $d \geq 2$. To the best of our knowledge, the only result in the literature in dimension greater than two is that of Marzo \cite{marzocubes}, who showed that there exists a Riesz basis of exponentials for $S$ whenever $S \subset \br^d$ is a finite union of disjoint translates of the unit cube. 

\subsection{}
In this paper we prove a result in spirit of Fuglede's theorem, but we consider the more general setting when $S$ tiles $\br^d$ \emph{with multiplicity}. The set $S$ is said to $k$-tile $\br^d$ with the translation set $\Lambda$ if almost
every point gets covered exactly $k$ times by the sets $S + \lambda$ $(\lambda \in \Lambda)$, that is,
\begin{equation}
\sum_{\lambda \in \Lambda} \1_S(x - \lambda) = k \quad \text{(a.e.)}
\label{eq:tilingconst}
\end{equation}
Clearly, the family of sets which $k$-tile $\br^d$ for some $k$ is much larger than the family of sets which just tile (i.e.\ $1$-tile). We will prove: 
\begin{theorem}
Let $S$ be a bounded, Riemann measurable set which $k$-tiles $\br^d$ with a lattice $\Lambda$. Then $S$ has a Riesz basis of exponentials.
\label{thm:multitiling}
\end{theorem}

By ``Riemann measurable'' we mean that the boundary of $S$ is a set of measure zero.

Let us mention some special cases of Theorem \ref{thm:multitiling} which may be of particular interest. The first one is a generalization of the above mentioned result due to Marzo \cite{marzocubes}.
\begin{corollary}
 Let $S$ be a bounded, Riemann measurable set which tiles $\br^d$ with a lattice $\Lambda$. Then any finite union of disjoint translates of $S$ has a Riesz basis of exponentials. 
\label{cor:unions}
\end{corollary}

The result in \cite{marzocubes} corresponds to the case when $S=[0,1)^d$ and $\Lambda = \bz^d$. Corollary \ref{cor:unions} shows that one may replace the cube with other fundamental domains of a lattice $\Lambda$. For example, any planar domain which can be tiled by finitely many translates of a hexagon has a Riesz basis of exponentials.

In dimension one, Corollary \ref{cor:unions} implies that a finite union of disjoint intervals with commensurable lengths has a Riesz basis of exponentials, a result proved in \cite{seip}. 

Theorem \ref{thm:multitiling} also allows us to extend the Lyubarskii-Rashkovskii result \cite{yurarash} on convex, symmetric polygons to dimensions greater than two. 
\begin{corollary}
 Let $S \subset \br^d$ be a centrally symmetric polytope, whose $(d-1)$-dimensional faces are also centrally symmetric, and whose vertices lie in some lattice $\Lambda$. Then $S$ admits a Riesz basis of exponentials.
\label{cor:symmpoly}
\end{corollary}

This is obtained by combining Theorem \ref{thm:multitiling} with a recent result \cite{gravin} that any centrally symmetric polytope in $\br^d$ whose $(d-1)$-dimensional faces are also centrally symmetric, and whose vertices lie in $\Lambda$, satisfies \eqref{eq:tilingconst} for some $k$.

We point out that in contrast to the result of \cite{yurarash} in two dimensions, in Corollary \ref{cor:symmpoly} we do not require that the polytope S is convex. On the other hand, the result in \cite{yurarash} is not completely covered by Corollary \ref{cor:symmpoly}, as we require that the vertices of the polytope lie in some lattice. 

Examples of polytopes satisfying the conditions in Corollary \ref{cor:symmpoly} include all zonotopes with vertices in a lattice $\Lambda$. Recall that a zonotope is the image of a cube in $\br^m$, $m \geq d$, under a linear mapping onto $\br^d$ (see \cite{ziegler}). Other examples in any dimension $d \geq 4$ may be found e.g.\ in \cite{mcmullen}. 

\section{Proof}

\subsection{}
To prove Theorem \ref{thm:multitiling} we use the approach due to Matei and Meyer \cite{meyer1, meyer3} based on simple quasicrystals. This approach was already used in \cite{kozmalev, lev} to construct a Riesz basis of exponentials on a union of intervals subject to an arithmetical condition on the lengths of the intervals.

Let $\Gamma$ be a lattice in $\br^{d+1} = \br^d \times \br$, and let $p_1$ and $p_2$ denote the projections onto $\br^d$ and $\br$, respectively. We assume that the restrictions of $p_1$ and $p_2$ to $\Gamma$ are injective, and that their images are dense. Denote by $\Gamma^*$ the dual lattice, consisting of all vectors $\gamma^* \in \br^{d+1}$ such that $\dotprod{\gamma}{\gamma^*} \in \bz$, $\gamma \in \Gamma$.

Let $S$ be a bounded, Riemann measurable set in $\br^d$, and let $I = [a, b)$ be a semi-closed interval. Define the Meyer ``cut-and-project'' sets
\[
\begin{aligned}
\Lambda(\Gamma, I) &= \{ p_1(\gamma) : \gamma \in \Gamma, \; p_2(\gamma) \in I\} \subset \br^d,\\[4pt]
\Lambda^*(\Gamma, S) &= \{ p_2(\gamma^*) : \gamma^* \in \Gamma^*, \; p_1(\gamma^*) \in S\} \subset \br .
\end{aligned}
\]
A key principle in \cite{meyer1, meyer3} is a ``duality'' between the quasicrystals $\Lambda(\Gamma, I)$ and $\Lambda^*(\Gamma, S)$, which allows us to reduce the problem to the single interval $I$.

\begin{lemma}
\label{lemma:duality}
Suppose that $E(\Lambda^*(\Gamma, S))$ is a Riesz basis in $L^2(I)$. Then $E(\Lambda(\Gamma, I))$ is a Riesz basis in $L^2(S)$.
\end{lemma}

The proof of Lemma \ref{lemma:duality} is along similar lines as in \cite{meyer3} (Sections 6--7). See also \cite{kozmalev} for a proof of the duality lemma in the periodic setting. 

\subsection{}
Now assume $S$ to $k$-tile $\br^d$ with some lattice. We may apply a linear transformation such that the $k$-tiling condition \eqref{eq:tilingconst} becomes
\begin{equation}
\sum_{m \in \bz^d} \1_S(x - m) = k \quad 
\label{eq:ztiling}
\end{equation}
almost everywhere. Remark that this implies that $\mes S=k$. 

Consider the lattice $\Gamma$ and its dual $\Gamma^*$ in $\br^d \times \br$ defined by
\begin{equation*}
\begin{aligned}
 \Gamma &= \{ ((\operatorname{Id}+\beta \alpha^{\top})m -\beta n, \, n-\alpha^{\top}m) : m \in \bz^d, n \in \bz \} \\
\Gamma^* &= \{ (m+\alpha n, \, (1+\beta^{\top}\alpha)n+\beta^{\top}m) : m \in \bz^d, n \in \bz \} ,
\end{aligned}
\end{equation*}
where $\alpha$ and $\beta$ are column vectors in $\br^d$ and $\operatorname{Id}$ denotes the $d \times d$ identity matrix. We choose the vector $\alpha = (\alpha_1, \alpha_2, \ldots , \alpha_d)^{\top}$ such that the numbers $1, \alpha_1, \alpha_2, \ldots , \alpha_d$ are linearly independent over the rationals, and the vector $\beta= (\beta_1, \beta_2, \ldots , \beta_d)^{\top}$ such that the numbers $\beta_1, \beta_2, \ldots, \beta_d, 1+\beta^{\top} \alpha$ are linearly independent over the rationals. These conditions ensure that the projections $p_1$ and $p_2$ restricted to $\Gamma$ are injective and that their images are dense.

Let $I=[0,k)$. We will prove that the exponential system $E(\Lambda(\Gamma, I))$ is a Riesz basis in $L^2(S)$. According to Lemma \ref{lemma:duality} it will be enough to check that $E(\Lambda^*(\Gamma, S))$ is a Riesz basis in $L^2(I)$. To prove this we will use the following result due to Avdonin.
\begin{theorem*}[Avdonin \cite{avdonin}]
Let $\{ \lambda_j  , \, j \in \bz \}$ be a sequence in $\br$ satisfying the following three conditions:
\begin{enumerate}[label=\textup{(\alph*)}]
\addtolength{\itemsep}{5pt}
\item
\label{item:separation}
$\{\lambda_j\}$ is a separated sequence, that is, $\inf_{i \neq j} |\lambda_i - \lambda_j| > 0$;
\item
\label{item:boundedness}
$\sup_{j} |\delta_j| < \infty$, where $\delta_j = \lambda_j - j/k$;
\item
\label{item:cancellation}
There is a constant $c$ and a positive integer $N$ such that
\begin{equation*}
\sup_{n \in \bz} \; \Big| \frac1{N} \sum_{j=n+1}^{n+N} \delta_{j} \; - \; c\Big| < \frac1{4k} \, .
\end{equation*}
\end{enumerate}
Then the system $\{ e^{2\pi i \lambda_j x} \}$ is a Riesz basis in $L^2[0,k)$.
\end{theorem*}

This is a generalization of Kadec's $1/4$ theorem which corresponds to the case when $N=1$. In fact, the above statement is a special case of the result proved in \cite{avdonin}.

\subsection{}
Define 
\begin{equation}
S_n = S \cap (n\alpha + \bz^d) , \quad \Lambda_n = \{ n + \dotprod{x}{\beta} \, : \, x \in S_n \} , \quad n \in \bz .  
\label{eq:setsn}
\end{equation}
Then $\{ \Lambda_n \}$ forms a partition of $\Lambda^*(\Gamma, S)$. By an appropriate translation of $S$ we may assume that \eqref{eq:ztiling} holds for all $x \in \bz \alpha$. Hence, each ``block'' $\Lambda_n$ contains exactly $k$ elements. We may thus choose an enumeration $\{\lambda_j, \, j \in \bz\}$ of the set $\Lambda^*(\Gamma, S)$ such that
\[
\Lambda_n = \{ \lambda_j \, : \, kn+1 \leq j \leq k(n+1)\}, \quad n \in \bz .
\]
We will show that for this enumeration, conditions \ref{item:separation}, \ref{item:boundedness} and \ref{item:cancellation} above are satisfied.

Condition \ref{item:separation} is a general property of cut-and-project sets. If $\Lambda^*(\Gamma, S)$ is not separated, there must exist a sequence $\{ \gamma^*_j \}\subset \Gamma^*$ such that $0 \neq p_2(\gamma^*_j) \to 0$ and $p_1(\gamma^*_j) \in S-S$. Since $S$ is bounded, infinitely many of the $\gamma^*_j$ must coincide with an element $\gamma^* \in \Gamma^*$ such that $p_2(\gamma^*) = 0$, a contradiction.  

To get condition \ref{item:boundedness}, for each $j$ choose $n = n(j)$ such that $\lambda_j \in \Lambda_n$. Then by \eqref{eq:setsn} we have $\lambda_j = n + \dotprod{x}{\beta}$ for some $x \in S_n$. Hence 
\begin{equation*}
|\delta_j| = |\lambda_j - j/k| = |\dotprod{x}{\beta}+(n-j/k)| \leq R \norm{\beta} + 1 , 
\end{equation*}
where $R=\sup \norm{x}$, $x \in S$. This confirms \ref{item:boundedness}. 

Finally we turn to establish \ref{item:cancellation}. We will show that there is a constant $c$ such that, for any $\varepsilon > 0$ and any $N > N(\varepsilon)$,  we have
\begin{equation}
\sup_{n \in \bz} \; \Big| \frac1{N} \sum_{j=kn+1}^{k(n+N)} \delta_{j} \; - \; ck \Big| < \varepsilon \, .
\label{eq:avdoninfit}
\end{equation}
Since the $\delta_j$ are already known to be bounded, this will imply \ref{item:cancellation}.

Consider the function
\[
\phi(x) = \sum_{m \in \bz^d} \dotprod{x+m}{\beta} \, \1_S(x+m) .
\]
This function is $1$-periodic, and hence may be viewed as a function on the torus $\bt^d = \br^d / \bz^d$. We have
\[
\sum_{j=kr+1}^{k(r+1)} \delta_j
= \sum_{j=kr+1}^{k(r+1)} (\lambda_j - r)
- \sum_{j=kr+1}^{k(r+1)} \Big(\frac{j}{k} - r\Big)
= \phi(r\alpha) - \frac{k+1}{2} .
\]
Taking the average for $r=n, n+1, \ldots , n+N-1$ yields 
\begin{equation*}
\frac{1}{N} \sum_{j=kn+1}^{k(n+N)} \delta_j = \frac{1}{N} \sum_{r=n}^{n+N-1} \phi(r\alpha) - \frac{k+1}{2}.
\end{equation*}
The numbers $1, \alpha_1, \alpha_2, \ldots , \alpha_d$ are linearly independent over the rationals, so the points $\{r \alpha \}$ are well-distributed on $\bt^d$. Since $\phi$ is a Riemann integrable function, we have 
\begin{equation*}
\sup_{n \in \bz} \left| \frac{1}{N} \sum_{r=n}^{n+N-1} \phi(r\alpha) - \int_{\bt^d} \phi(x) dx \right| = o(1), \quad N \to \infty 
\end{equation*}
(see \cite[pp. 46, 52]{kuipers}). It follows that \eqref{eq:avdoninfit} holds with 
\begin{equation*}
c = \frac{1}{k} \int_{\bt^d} \phi(x) dx - \frac{k+1}{2k}.
\end{equation*}
This implies \ref{item:cancellation}, and concludes the proof of Theorem \ref{thm:multitiling}.

\section{Remarks}
1.\ We point out that the Riesz basis constructed for the multi-tiling set $S$ can be seen to depend only on the translation lattice $\Lambda$ and the tiling multiplicity $k$ in condition \eqref{eq:tilingconst}, and not on the specific structure of the set $S$. 

2.\ In the one-dimensional case $S \subset \br$, Theorem \ref{thm:multitiling} can also be deduced from Lemma 2 in \cite{kozmanitzan}.

3.\ Shortly after the preliminary version of this paper was posted, a simpler proof of Theorem \ref{thm:multitiling} was found by Kolountzakis \cite{koloun1}. 

\subsection*{Acknowledgement}
We would like to thank the Universitat de Barcelona and the Centre de Recerca Matem\`atica (CRM) for hosting and supporting us during the research stay that led to this collaboration. 


\end{document}